\documentclass[leqno,12pt]{amsart}

\usepackage{amssymb}
\usepackage{amscd}              
\usepackage[all]{xy}            
\xyoption{curve}
\CompileMatrices

\newcommand{\Ac}{\mathcal{A}}

\newcommand{\bO}{\text{\bf O}}

\newcommand{\C}{\mathbb{C}}

\newcommand{\card}{\sharp}
\newcommand{\chara}{\text{char}}                

\newcommand{\comp}{\circ}
\newcommand{\congruent}{\equiv}
\newcommand{\Ct}{\widetilde{C}}                  

\newcommand{\diag}{\text{diag}}

\newcommand{\eins}{{\mathbb{I}}}

\newcommand{\E}{{\mathcal{E}}}

\newcommand{\etah}{\widehat{\eta}}

\newcommand{\F}{{\mathcal{F}}}

\newcommand{\Ft}{\tilde{\F}}

\newcommand{\gen}{\text{gen}}

\newcommand{\GL}{\text{GL}}

\newcommand{\GVB}{\text{GVB}}

\newcommand{\id}{\text{id}}

\newcommand{\injto}{\hookrightarrow}
\newcommand{\isomorph}{\cong}           
\newcommand{\isomto}{\overset{\sim}{\rightarrow}}

\newcommand{\Isom}{\text{Isom}}
\newcommand{\Isomto}{\overset{\sim}{\longrightarrow}}

\newcommand{\KGL}{\text{KGL}}

\newcommand{\Kh}{\widehat{K}}

\newcommand{\m}{{\mathfrak{m}}}                  
\newcommand{\M}{{\mathcal{M}}}

\newcommand{\Oh}{\widehat{\Oo}}
\newcommand{\Oo}{{\mathcal{O}}}                  

\newcommand{\PB}{\text{PB}}

\newcommand{\Pp}{{\mathbb{P}}}                   

\newcommand{\psib}{\overline{\psi}}

\newcommand{\sh}{\text{sh}}
\newcommand{\sing}{\text{sing}}
\newcommand{\Spec}{\text{Spec}\, }               

\newcommand{\tensor}{\otimes}
\newcommand{\TVB}{\text{TVB}}

\newcommand{\tF}{\tilde{F}}

\newcommand{\tH}{\tilde{H}}

\newcommand{\V}{{\mathcal V}} 
\newcommand{\varphib}{\overline{\varphi}}

\newcommand{\Z}{{\mathbb{Z}}}                    

\newtheorem{theorem}{Theorem}[section]
\newtheorem{proposition}[theorem]{Proposition}
\newtheorem{lemma}[theorem]{Lemma}
\newtheorem{corollary}[theorem]{Corollary}

\theoremstyle{definition}
\newtheorem{definition}[theorem]{Definition}

\newenvironment{inmargins}[1]{\begin{list}{}{
    \leftmargin=#1 \rightmargin=0pt \parsep=0pt
    \partopsep=0pt}\item[]}{\end{list}}

\begin{document}

\title[Twisted vector bundles on pointed nodal curves]
{Twisted vector bundles on pointed nodal curves}
\author[Ivan Kausz]{Ivan Kausz}
\thanks{Partially supported by the DFG} 
\date{January 12, 2005}
\address{NWF I - Mathematik, Universit\"{a}t Regensburg, 93040 Regensburg, 
Germany}
\email{ivan.kausz@mathematik.uni-regensburg.de}
\begin{abstract}
Motivated by the quest for a good compactification
of the moduli space of $G$-bundles on a nodal curve
we establish a striking relationship between Abramovich's and
          Vistoli's twisted bundles and Gieseker vector bundles.
\end{abstract}


\maketitle

\tableofcontents


\section{Introduction}

This paper grew out of an attempt to understand a recent 
draft of Seshadri (\cite{Seshadri-notes})
and is meant as a contribution in the quest for a good compactification
of the moduli space (or stack) of $G$-bundles on a nodal curve.

We are led by the idea that such a compactification should behave well in
families and also under partial normalization of nodal curves.
This statement may be reformulated by saying that we are looking for
an object which has the right to be called
the moduli stack of stable maps into the classifying stack $BG$
of a reductive group $G$.

For finite groups $G$ the stack of stable maps into $BG$ has been 
recently constructed by means of so called twisted bundles
by D. Abramovich and A. Vistoli
(\cite{AV}, \cite{ACV}).  
On the other hand, as shown in
\cite{stable maps}, the notion of Gieseker vector bundles leads to 
the construction of the stack of stable maps into $B\GL_r$. 

In this note we establish a connection between the straightforward
generalization of the notion of twisted bundles to the case of the non-finite
reductive group $\GL_r$ and Gieseker vector bundles.
My hope is that this relationship - whose observation is entirely due
to Seshadri, and which in my mind is really striking -
may help to find the right notion for more general reductive
groups $G$.

I would like to thank Seshadri for generously imparting his ideas.
This paper owes very much to long discussions which I had with
Nagaraj in November and December 2002. I would like to thank the Institute
of Mathematical Sciences in Chennai, whose hospitality made these
discussions possible.

\vspace{5mm}
\section{Twisted $G$-bundles}

Throughout this section $k$ denotes an algebraically closed field
and $G$ a reductive group over $k$.

A twisted $G$-bundle is a twisted object in the sense of 
\cite{AV}, \S 3 where the target stack $\M$ is taken 
to be the classifying stack $BG$.
For convenience we recall the necessary definitions from loc. cit.

\begin{definition}
1.
An {\em $n$-marked curve} consists of data $(U\to S,\Sigma_i)$
where $\pi:U\to S$ is a nodal curve and 
$\Sigma_1,\dots,\Sigma_n\subset U$ are pairwise disjoint closed
subschemes whose supports do not intersect the singular locus
$U_{\sing}$ of $\pi$ and are such that the projections $\Sigma_i\to S$ are
\'etale.

2. 
A {\em morphism between two $n$-marked curves 
$(U\to S, \Sigma^U_i)$ and $(V\to S, \Sigma^V_i)$}
is an $S$-morphism $f:U\to V$ such that 
$f(\Sigma^U_i)\subseteq \Sigma^V_i$ for each $i$.
Such a morphism is called {\em strict}, if for each $i$ the support of 
$f^{-1}(\Sigma^V_i)$  coincides with the support of $\Sigma^U_i$
and if furthermore 
the support of $f^{-1}(V_{\sing})$ coincides with the one of $U_{\sing}$.

3.
The {\em pull back} of an $n$-marked curve
$(U\to S, \Sigma_i)$ by a morphism $S'\to S$ is the $n$-marked curve
$(U\times_SS', \Sigma_i\times_SS')$.

4.
An {\em $n$-pointed nodal curve} is an $n$-marked curve where
the projections $\Sigma_i\to S$ are isomorphisms.

5.
Let $(U\to S,\Sigma_i)$ be an $n$-marked curve.
The complement (inside $U$) of the union of the singular
locus $U_{\sing}$ and the markings $\Sigma_i$  is called 
the {\em generic locus} of $U$ and is denoted by $U_{\gen}$.
\end{definition}

\begin{definition}
1.
An {\em action of a finite group $\Gamma$} on an
$n$-marked nodal curve $(U\to S, \Sigma_i)$  is an action
of $\Gamma$ on $U$ as an $S$-scheme which leaves the $\Sigma_i$ invariant.
Such an action is called {\em tame}, if for each geometric point
$u$ of $U$ the stabilizer $\Gamma_u\subseteq\Gamma$ of $u$ has order prime to
the characteristic of $u$.

2. 
Let $S$ be a $k$-scheme.
Let $(U\to S, \Sigma_i)$ be an $n$-marked nodal curve and let 
$\eta$ be a principal $G$-bundle on $U$.
A {\em essential action of a finite group $\Gamma$ on $(\eta, U)$}  
is a pair of actions of $\Gamma$ on $\eta$ and on $(U\to S, \Sigma_i)$ 
such that
\begin{enumerate}
\def\theenumi{\roman{enumi}}
\def\labelenumi{(\theenumi)}
\item
the actions of $\Gamma$ on $\eta$ and on $U$ are compatible, i. e.
if $\pi: \eta\to U$ denotes the projection, then 
$\pi\comp\gamma = \gamma\comp\pi$ for each $\gamma\in\Gamma$.
\item
if $\gamma\in\Gamma$ is an element different from the identity and
$u$ is a geometric point of $U$ fixed by $\gamma$, then
the automorphism of the fiber $\eta_u$ induced by $\gamma$ is not trivial.
\end{enumerate}

3.
An essential action of a finite group $\Gamma$ on $(\eta, U)$ is
called {\em tame}, if the action of $\Gamma$ on $(U\to S, \Sigma_i)$ 
is tame.
\end{definition}

\begin{definition}
Let $S$ be a $k$-scheme.
Let $C\to S$ be an $n$-pointed nodal curve and let $\xi$ be a principal
$G$-bundle over $C_{\gen}$. 
A {\em chart  $(U,\eta,\Gamma)$ for $\xi$} consists of the following
data
\begin{enumerate}
\item
An $n$-marked curve $U\to S$ and a strict morphism $\phi: U\to C$,
\item
A principal $G$-bundle $\eta$ on $U$.
\item
An isomorphism $\eta\times_UU_{\gen}\isomto \xi\times_CU_{\gen}$
of $G$-bundles on $U_{\gen}$.
\item
A finite group $\Gamma$.
\item
A tame, essential action of $\Gamma$ on $(\eta, U)$.
\end{enumerate}

These data are required to satisfy the following conditions

\begin{enumerate}
\def\theenumi{\roman{enumi}}
\def\labelenumi{(\theenumi)}
\item
The action of $\Gamma$ leaves the morphisms $U\to C$ and
$\eta\times_UU_{\gen}\isomto \xi\times_CU_{\gen}$ invariant.
\item
The induced morphism $U/\Gamma\to C$ is \'etale.
\end{enumerate}
\end{definition}

\begin{proposition}
\label{cyclic}
(cf. \cite{AV}, Prop 3.2.3)
Let $C\to S$ be an $n$-pointed nodal curve over a $k$-scheme $S$
and let $\xi$ be a principal $G$-bundle on $C_{\gen}$.
Let $(U,\eta,\Gamma)$ be a chart for $\xi$. Then the following
holds.
\begin{enumerate}
\item
The action of $\Gamma$ on $U_{\gen}$ is {\em free}. 
\end{enumerate}
Let $s$ be a geometric point of $S$ and let $u$ be a closed
point of the curve $U_s$. Let $\Gamma_u\subseteq \Gamma$ be the 
stabilizer of $u$. Then $\Gamma_u$ is a cyclic group.
Let $e$ be its order and let $\gamma_u$ be a generator of $\Gamma_u$.
Then
\begin{enumerate}
\setcounter{enumi}{1}
\item
if $u$ is a regular point, the action of $\gamma_u$ on the tangent
space of $U_s$ at $u$ is via multiplication by a primitive $e$-th 
root of unity.
\item
if $u$ is a singular point, $\Gamma_u$ 
leaves each of the two branches of $U_s$ at $u$ invariant.
The action of $\gamma_u$ on the tangent space of each of the branches is
via multiplication with a primitive $e$-th root of unity. 
\end{enumerate}
\end{proposition}


\begin{definition}
Let $C\to S$ be an $n$-pointed nodal curve over a $k$-scheme $S$ and
let $\xi$ be a principal $G$-bundle on $C_{\gen}$.
A chart $(U,\eta,\Gamma)$ for $\xi$ is called {\em balanced},
if for each geometric fiber of $U\to S$ and each singular point $u$ on 
it the action of $\gamma_u$ on the tangent spaces of the two branches
is via multiplication with primitive roots of unity which are inverse 
to each other.
\end{definition}

\begin{definition}
Let $C\to S$ be an $n$-pointed nodal curve over a $k$-scheme $S$ and
let $\xi$ be a principal $G$-bundle on $C_{\gen}$.
Two charts $(U_1, \eta_1, \Gamma_1)$ 
and $(U_2, \eta_2, \Gamma_2)$ of $\xi$ are called {\em compatible},
if for each pair of $u_1$, $u_2$ of geometric points of $U_1$, $U_2$
lying above the same geometric point $u$ of $C$ the following holds:
\begin{inmargins}{8mm}
Let $C^{\sh}$ denote the strict henselization of $C$ at $u$.
For $j=1,2$ 
let $\Gamma'_j\subseteq\Gamma_j$ denote the stabilizer
subgroup of the point $u_j$, 
let $U_j^{\sh}$ denote the strict henselization
of $U_j$ at $u_j$, and let $\eta_j^{\sh}:=\eta_j\times_{U_j}U_j^{\sh}$.
Then there exists an isomorphism $\theta: \Gamma'_1\to\Gamma'_2$,
a $\theta$-equivariant isomorphism $\phi:U_1^{\sh}\isomto U_2^{\sh}$
of $C^{\sh}$-schemes, and a $\theta$-equivariant isomorphism
$\eta_1^{\sh}\isomto\phi^*\eta_2^{\sh}$ of $G$-bundles.
\end{inmargins}
\end{definition}

\begin{definition}
Let $g$ and $n$ be two non-negative integers.
An  {\em $n$-pointed twisted $G$-bundle of genus $g$} is a triple
$(\xi, C\to S, \Ac)$ where

\begin{enumerate}
\item
$S$ is a $k$-scheme,
\item
$C\to S$ is proper $n$-pointed nodal curve of finite presentation
with geometrically connected fibers of genus $g$,
\item
$\xi$ is a principal $G$-bundle on $C_{\gen}$,
\item
$\Ac=\{(U_\alpha,\eta_\alpha,\Gamma_\alpha)\}$ is a balanced atlas, i.e.
a collection of
mutually compatible balanced charts for $\xi$, such that the images of the
$U_\alpha$ cover $C$.
\end{enumerate}

\end{definition}

\begin{definition}
Let $(\xi, C\to S, \Ac)$ be an $n$-pointed twisted $G$-bundle
of genus $g$. A morphism of $k$-schemes $S'\to S$ induces
a triple $(\xi', C'\to S', \Ac')$ as follows:

\begin{itemize}
\item
The $n$-pointed nodal curve $C'\to S'$ is the pull back of 
$C\to S$ by $S'\to S$.
\item
Thus we have a morphism $C'_{\gen}\to C_{\gen}$, and the $G$-bundle 
$\xi'$ is the pull back of $\xi$ by this morphism.
\item
Let $\{U_\alpha,\eta_\alpha,\Gamma_\alpha)\}$ be the set of charts 
which make up the atlas $\Ac$. Then 
$\Ac'=\{U'_\alpha,\eta'_\alpha,\Gamma_\alpha)\}$, where
$U_{\alpha}'\to S'$ is the pull back of the $n$-marked curve 
$U_{\alpha}\to S$,
and $\eta'_\alpha$ is the pull back of $\eta_\alpha$ by the morphism
$U'_{\alpha}\to U_{\alpha}$.
Since the $(U'_\alpha, \eta'_\alpha, \Gamma_\alpha)$ are charts for
$\xi'$ which are balanced and mutually compatible (cf. \cite{AV},
Prop. 3.4.3), $\Ac'$ is a balanced atlas.
\end{itemize}

Thus the triple $(\xi', C'\to S', \Ac')$ is an
$n$-pointed twisted $G$-bundle of genus $g$. It is called the
{\em pull back} of $(\xi, C\to S, \Ac)$ by the morphism $S'\to S$.
\end{definition}

\begin{definition}
A {\em morphism} between two $n$-pointed twisted $G$-bundles 
$(\xi', C'\to S', \Ac')$ and $(\xi, C\to S, \Ac)$ consists
of a Cartesian diagram
$$
\xymatrix@R=2ex{
C' \ar[r] \ar[d] &
C \ar[d] \\
S' \ar[r] &
S}
$$
and an isomorphism $\xi'\isomto\xi\times_{C_{\gen}}C'_{\gen}$
such that the pull-back of the charts in $\Ac$ 
(considered as charts for $\xi'$) are compatible
with all the charts in $\Ac'$.
\end{definition}

\vspace{5mm}
\section{Review of Gieseker vector bundles}
\label{review}

In this section I will recall some definitions from my earlier papers
\cite{kgl} and \cite{degeneration}.

Let $k$ be an algebraically closed field. 
Let $n\geq 1$ be an integer and let $R_1,\dots,R_n$ be $n$
copies of the projective line $\Pp^1$. On each $R_i$ we choose
two distinct points $x_i$ and $y_i$.
Let $R$ be the nodal curve over $k$ constructed from $R_1,\dots, R_n$
by identifying $y_i$ with $x_{i+1}$ for $i=1,\dots,n-1$.
We call such a curve $R$ a {\em chain of projective lines} of length $n$ with components
$R_1,\dots,R_n$. On the extremal components $R_1$ and $R_n$ we have the
two points $x_1$ and $y_n$ respectively, which are smooth points of $R$.

\begin{definition}
\label{admissible}
A vector bundle $E$ of rank $r$ on $R$ is called {\em admissible},
if 
\begin{enumerate}
\item
for each $i\in[1,n]$ the restriction of $E$ on the component $R_i$ is
of the form
$$
d_i\Oo_{R_i}(1)\oplus(r-d_i)\Oo_{R_i}
$$
for some integer $d_i\geq 1$ and
\item
there exists no nonvanishing global section of $E$ over $R$ which
vanishes in the two points $x_1$ and $y_n$.
\end{enumerate}
\end{definition}

Let $C$ be an irreducible curve with exactly one double point $p$.
Let $\Ct\to C$ be the normalization of $C$ and let $p_1, p_2\in \Ct$
be the two points lying above $p$.
Let $C_0:=C$. 
For $n\geq 1$ we let $C_n$ denote reducible nodal curve
which is constructed from $\Ct$ and a chain $R=R_1\cup\dots\cup R_n$
of projective lines by identifying the points $p_1, x_1$ and
$p_2, y_n$ respectively.

\begin{definition}
\label{GVB}
A {\em Gieseker vector bundle} on $C$ is a pair $(C'\to C,\F)$ where
$C'=C_n$ for some $n\geq 0$, the morphism $C'\to C$ is the one which
contracts the chain of projective lines to the point $p$ and $\F$ is
a vector bundle on $C'$ whose restriction to the chain of projective
lines is admissible in the sense of \ref{admissible}.
\end{definition}

\begin{definition}
\label{GVBD}
A {\em Gieseker vector bundle datum} on the two-pointed curve
$(\Ct,p_1,p_2)$ is a triple $(C'\to C, F, p')$, where 
$(C'\to C,\F)$ is a Gieseker vector bundle on $C$ and $p'$ is
a singular point in $C'$. 
\end{definition}

Let $V$ and $W$ be two $r$-dimensional $k$-vector spaces.
In \cite{kgl} I have constructed a certain compactification
$\KGL(V,W)$
of the space $\Isom(V,W)$ of linear isomorphisms from $V$ to $W$
which has properties similar to De Concinis and Procesis so called
wonderfull compactification of adjoint linear groups.
We need the following fact about $\KGL(V,W)$ whose proof can
be found in \cite{kgl}, \S 9:

The variety $\KGL(V,W)$
is the disjoint union of strata $\bO_{I,J}\subset\KGL(V,W)$ indexed
by pairs of subsets $I,J\in[0,r-1]$ such that $\min(I)+\min(J)\geq r$.
Let $I,J$ be such a pair.
Let us write $I=\{i_1,\dots,i_{n_1}\}$ and $J=\{j_1,\dots,j_{n_2}\}$
where $i_1<\dots<i_{n_1}<i_{n_1+1}:=r$ and $j_1<\dots<j_{n_2}<j_{n_2+1}:=r$.
A ($k$-valued) point in $\bO_{I,J}$ is given by
the data 
$$
\Phi=
(F_{\bullet}(V),F_{\bullet}(W),\varphib_1,\dots,\varphib_{n_1},\psib_1,\dots,\psib_{n_2},\Phi')
$$ 
where
\begin{enumerate}
\item
$F_{\bullet}(V)$ denotes a flag 
$$
0=F_0(V)\subsetneq F_1(V)\subsetneq\dots\subsetneq F_{n_2}(V)\subseteq F_{n_2+1}(V)
        \subsetneq\dots\subsetneq F_{n_1+n_2+1}(V)=V
$$
whith
$\dim F_{\nu}(V)=r-j_{n_2+1-\nu}$ for $\nu\in [0,n_2]$ and
$\dim F_{\nu}(V)=i_{\nu-n_2}$ for $\nu\in [n_2+1,n_1+n_2+1]$,
\item
$F_{\bullet}(W)$ denotes a flag 
$$
0=F_0(W)\subsetneq F_1(W)\subsetneq\dots\subsetneq F_{n_1}(W)\subseteq F_{n_1+1}(W)
        \subsetneq\dots\subsetneq F_{n_1+n_2+1}(W)=W
$$
where 
$\dim F_{\nu}(W)=r-i_{n_1+1-\nu}$ for $\nu\in [0,n_1]$ and
$\dim F_{\nu}(W)=i_{\nu-l}$ for $\nu\in [n_1+1,n_1+n_2+1]$,
\item
the symbol
$\varphib_{\nu}$ denotes the homothety class of an isomorphism
from  the subquotient $F_{n_1-\nu+1}(W)/F_{n_1-\nu}(W)$ of $W$ to the subquotient
$F_{n_2+\nu+1}(V)/F_{n_2+\nu}(V)$ of $V$,
\item
the symbol
$\psib_{\nu}$ denotes the homothety class of an isomorphism
from  the subquotient $F_{n_2-\nu+1}(V)/F_{n_2-\nu}(V)$ of $V$ to the subquotient
$F_{n_1+\nu+1}(W)/F_{n_1+\nu}(W)$ of $W$,
\item
the symbol 
$\Phi'$ denotes an isomorphism from the subquotient
$F_{n_2+1}(V)/F_{n_2}(V)$ of $V$ to the subquotient $F_{n_1+1}(W)/F_{n_1}(W)$ of $W$.
\end{enumerate}

The relationship between Gieseker vector bundles and the compactification $KGL(V,W)$ is
given by the following

\begin{theorem}
\label{GVBD->KGL}
(Cf. \cite{degeneration}, Theorem 9.5) 
There exists a natural
bijection from the set of all Gieseker vector bundle data on $(\Ct,p_1,p_2)$ 
to the set of all pairs $(\E,\Phi)$, where $\E$ is a vector bundle on
$\Ct$ and $\Phi$  is a $k$-valued point in 
$
\KGL(\E[p_1],\E[p_2])
$.

More precisely,
let  $(C'\to C,\F,p')$ be a Gieseker vector bundle datum on $(\Ct,p_1,p_2)$.
Let $R=R_1\cup\dots\cup R_n$ be the chain of projective lines in $C'$.
Let $y_0:=p_1$ and $x_{n+1}:=p_2$. Let $n_1+n_2=n$ be such that
the singular point $p'\in C'$ comes from identifying the points $y_{n_1}$ and
$x_{n_1+1}$.
Let $d_i$ be the degree of $\F$ restricted to $R_i$. 
Let $(\E,\Phi)$ be the pair associated to the given Gieseker vector bundle datum
$(C'\to C, \F, p')$.
Then $\Phi$ is in fact a point in the stratum
$\bO_{I,J}$, where $I=\{i_1,\dots,i_{n_1}\}$, $J=\{j_1,\dots,j_{n_2}\}$ and the 
$i_{\nu}$, $j_{\nu}$ are defined by
$$
i_{\nu} = r-\sum_{i=\nu}^{n_1}d_i
\qquad,\qquad
j_{\nu} = r-\sum_{i=n_1+1}^{n-\nu+1}d_i
\quad.
$$
The special case $n=0$ is included here in the sense that then $I=J=\emptyset$
and $\Phi\in\bO_{\emptyset,\emptyset}=\Isom(\E[p_1],\E[p_2])$.
\end{theorem}

\vspace{5mm}
\section{Twisted $\GL_r$-bundles on a fixed curve}

Throughout this section $k$ denotes an algebraically closed field
and $r$ a positive integer.

Let $(C,p_i)$ be an $n$-pointed nodal curve over $k$. 
Let $\TVB_r(C,p_i)$ be the set of isomorphism classes of
$n$-pointed twisted $\GL_r$-bundles of the form
$$
(\xi, C\to\Spec(k), \Ac)
\quad.
$$

\vspace{3mm}
\noindent
{\bf The case of a one-pointed smooth curve.}
Assume that $C$ is smooth and that $n=1$, i.e. $(C,p_i)=(C,p)$ is a
one-pointed smooth curve. Let $\PB_r(C,p)$ be the set of isomorphism
classes of vector bundles $E$ of rank $r$ on $C$ together with a
flag in the fiber at $p$.

\begin{theorem}
\label{thm1}
There is a natural surjection
$$
\TVB_r(C,p)\to\PB_r(C,p)
\quad.
$$
\end{theorem}

We skip the proof of Theorem \ref{thm1}, since on the one hand the result 
is well known (cf. \cite{MS}, \cite{Biswas}) and on the other hand 
there is a proof analogous to (and easier than) the proof of Theorem 
\ref{thm2} which we give in detail below.

\vspace{3mm}
\noindent
{\bf The case of a nodal curve with one singularity.}
Assume now that $n=0$ and $C$ has exactly one double point.
Let $\GVB_r(C)$ be the set of isomorphism classes of
Gieseker vector bundles of rank $r$ on $C$.

\begin{theorem}
\label{thm2}
There is a natural surjection
$$
\TVB_r(C)\to\GVB_r(C)
\quad.
$$ 
\end{theorem}

The rest of the paper is concerned with the proof of Theorem \ref{thm2}.

\vspace{5mm}
\section{Construction}

Let $C$ be a nodal curve over $\Spec(k)$ with one singular point $p$.
Let $(\xi, C\to\Spec(k), \Ac)$ be an object of $\TVB_r(C)$.
Let $(U,\eta,\Gamma)$ be a chart belonging to $\Ac$ such that
there is a point $q\in U$ which is mapped to $p$.

We denote by $\widehat{\Oo}_{p}$ and $\widehat{\Oo}_{q}$ 
the completion of the local rings $\Oo_{C,p}$ and $\Oo_{U,q}$
respectively. 
Let $\Gamma_q\subseteq\Gamma$ be the subgroup consisting of those
elements, which leave $q$ invariant. 
$\Gamma_q$ acts on $\widehat{\Oo}_q$, and $\widehat{\Oo}_p$
may be identified with the set of invariants under that action.
By proposition \ref{cyclic} the group $\Gamma_q$ is cyclic of
some order $e$ (which is prime to $\chara(k)$ by the tameness
assumption). 
Let $\gamma$ be a generator of $\Gamma_q$.

We choose an isomorphism
\begin{eqnarray} 
\label{1}
\widehat{\Oo}_p\isomto k[[s,t]]/(s\cdot t)
\quad.
\end{eqnarray}

It follows from \ref{cyclic}.(3) that there exists an isomorphism
\begin{eqnarray}
\label{2}
\widehat{\Oo}_q\isomto k[[u,v]]/(u\cdot v)
\end{eqnarray}
and a primitive $e$-th root of unity $\zeta$ such that the
diagrams
$$
\xymatrix{
\text{$\widehat{\Oo}_q$} \ar[r]^(0.3){\isomorph} 
                         \ar[d]^{\gamma}
&
k[[u,v]]/(u\cdot v) \ar[d] &
u \ar@{|->}[d]  & v \ar@{|->}[d] 
\\
\text{$\widehat{\Oo}_q$} \ar[r]^(0.3){\isomorph} 
&
k[[u,v]]/(u\cdot v) &
\zeta u & \zeta^{-1}v 
}
$$
and 
$$
\xymatrix{
\text{$\widehat{\Oo}_q$} \ar[r]^(0.3){\isomorph} 
&
k[[u,v]]/(u\cdot v) &
u^e & v^e
\\
\text{$\widehat{\Oo}_p$} \ar[r]^(0.3){\isomorph} 
                         \ar@{^(->}[u]
&
k[[s,t]]/(s\cdot t) \ar@{^(->}[u] &
s \ar@{|->}[u] & t \ar@{|->}[u] 
}
$$
are commutative.

Let $\widehat{K}_p$ be the total quotient ring of $\widehat{\Oo}_p$.
Then we have
$
\Spec(\widehat{K}_p)=\Spec(\widehat{\Oo}_p)\times_CC_{\gen}
$ 
and the isomorphism (\ref{1}) induces an isomorphism
$
\widehat{K}_p \isomto k((s))\times k((t))
$.
We choose an isomorphism
\begin{eqnarray}
\label{3}
\xi\times_{C_{\gen}}\Spec(\widehat{K}_p)\isomto
\GL_r\times\Spec(\widehat{K}_p)
\quad.
\end{eqnarray}

The group $\Gamma_q$ acts on $\eta\times_U\Spec(\widehat{\Oo}_q)$
(since it acts compatibly on $\eta$, $U$, $\Spec(\widehat{\Oo}_q)$).
To analyse this action we need the following

\begin{lemma}
\label{action}
Let $k$ be an algebraically closed field.
Let $(R,\m)$ be a local $k$-algebra with residue field $R/\m=k$.
Let $\Gamma$ be a cyclic group of order $e$ prime to the
characteristic of $k$ and let $\gamma\in\Gamma$ be a generator.
Assume that $\Gamma$ acts on $R$ such that the induced action
on $k$ is trivial.
Let $M$ be a trivial $R$-module of rank $r$ on which $\Gamma$
acts such that $\gamma(ax)=\gamma(a)\gamma(x)$ for all 
$a\in R$, $x\in M$.
Then there is a basis $x_1,\dots,x_r$ of $M$ such that
$\gamma(x_i)=\zeta_ix_i$ for some $e$-th roots of unity $\zeta_i$.
\end{lemma}

\begin{proof}
Let $e_1,\dots,e_r$ be an arbitrary basis of $M$.
Let $a=(a_{i,j})\in\GL_r(R)$ be defined by 
$
\gamma(e_j)=\sum_{i=1}^ra_{i,j}e_i
\quad.
$
Since $\gamma$ is of order $e$, it follows that
$$
\prod_{j=0}^{e-1}\gamma^j(a)=1
\quad.
$$
We have to show that there is a matrix $b\in\GL_r(R)$ 
such that
$$
a\cdot\gamma(b)=b\cdot z
$$
for some diagonal matrix $z\in\GL_r(k)$ with $z^e=1$.

Representation theory of finite groups tells us
that there is a matrix $c\in\GL_r(k)$ and a diagonal matrix
$z\in\GL_r(k)$ with $z^e=1$ such that 
$a\cdot c\equiv c\cdot z$ modulo $\m$.
Let $a':=c^{-1}\cdot a\cdot c$ and let $b'$ be the matrix
$$
b':=\sum_{i=0}^{e-1}\left(\prod_{j=0}^{i-1}\gamma^j(a')\right)z^{-i}
\quad.
$$
Since $b\equiv e\cdot 1$ modulo $\m$, it follows that $b'\in\GL_r(R)$.
Using the fact that $\prod_{i=0}^{e-1}\gamma^i(a')=1$ a simple
calculation shows that
$$
\gamma(b')=(a')^{-1}\cdot b'\cdot z
\quad.
$$
Therefore, if we set $b:=c\cdot b'$, we get the desired equality.
\end{proof}

\begin{corollary}
\label{alpha} 
There exists an isomorphism 
\begin{eqnarray}
\label{4}
\eta\times_U\Spec(\widehat{\Oo}_q)\isomto
\GL_r\times\Spec(\widehat{\Oo}_q)
\end{eqnarray}
of principal $\GL_r$-bundles on $\Spec(\widehat{\Oo}_q)$,
and elements $\alpha_1,\dots\alpha_r\in\Z/e\Z$
such that the following diagram commutes: 
$$
\xymatrix{
\text{$\eta\times_U\Spec(\widehat{\Oo}_q)$} 
\ar[r]^{\isomorph}
\ar[d]^{\gamma}
&
\GL_r\times\Spec(\widehat{\Oo}_q)
\ar[d]^{\diag(\zeta^{\alpha_1},\dots,\zeta^{\alpha_r})\times\gamma}
\\
\text{$\eta\times_U\Spec(\widehat{\Oo}_q)$} 
\ar[r]^{\isomorph}
&
\GL_r\times\Spec(\widehat{\Oo}_q)
}
$$
where the morphism 
$\diag(\zeta^{\alpha_1},\dots,\zeta^{\alpha_r}): \GL_r\to \GL_r$ 
is multiplication from the left with the matrix whose only
non-zero entries are the values $\zeta^{\alpha_1},\dots,\zeta^{\alpha_r}$
on the diagonal.
\end{corollary}

\begin{proof}
This is immediate from lemma \ref{action}.
\end{proof}

Let $\widehat{K}_q$ be the total quotient ring of 
$\widehat{\Oo}_q$.
The $\Gamma$-equivariant isomorphism 
$
\eta\times_UU_{\gen} \isomto \xi\times_{C_{\gen}}U_{\gen}
$,
which is part of the data of the chart $(U,\eta,\Gamma)$,
induces a $\Gamma_q$-equivariant isomorphism
\begin{eqnarray}
\label{5}
\eta\times_U\Spec(\widehat{K}_q)
\isomto
\xi\times_{C_{\gen}}\Spec(\widehat{K}_q)
\end{eqnarray}
of principal $\GL_r$-bundles over $\Spec(\widehat{K}_q)$.

Via the isomorphisms (\ref{3}) and (\ref{4}) such an isomorphism
is given by a matrix 
$
F \in \GL_r(\widehat{K}_q)
$
such that 
$$
\gamma(F) =
F \cdot \diag(\zeta^{\alpha_1},\dots,\zeta^{\alpha_r}) 
$$

The isomorphism (\ref{2}) induces an isomorphism 
$
\GL_r(\widehat{K}_q)\isomto \GL_r(k((u)))\times\GL_r(k((v)))
$
and we denote by $(F^1(u), F^2(v))$ the image of $F$ under this isomorphism.
The above condition on $F$ translates into the condition
\begin{eqnarray}
\label{6}
F^1_{i,j}(\zeta u) 
&=& 
\zeta^{\alpha_j} F^1_{i,j}(u) 
\\
\label{7}
F^2_{i,j}(\zeta^{-1} v) 
&=& 
\zeta^{\alpha_j} F^2_{i,j}(v)
\end{eqnarray} 
for the entries 
$F^1_{i,j}(u)\in k((u))$ and $F^2_{i,j}(v)\in k((v))$ 
of the matrices $F^1(u)$ and $F^2(v)$.

After possibly changing the isomorphism (\ref{4}) by a permutation
matrix, we can choose integers $a_1,\dots,a_r$ with
\begin{eqnarray}
\label{8}
0\leq a_1\leq a_2\leq \dots \leq a_r < e
\quad
\text{and $a_i\congruent\alpha_i$ mod $e\Z$.}
\end{eqnarray}

Conditions (\ref{6}), (\ref{7}) imply that there are matrices
$H^1(s)$ and $H^2(t)$ with entries 
$H^1_{i,j}(s)\in k((s))$ and $H^2_{i,j}(t)\in k((t))$
such that
\begin{eqnarray}
\label{9}
F^1_{i,j}(u) &=& u^{a_j} H^1_{i,j}(u^e)
\\
\label{10}
F^2_{i,j}(v) &=& v^{-a_j} H^2_{i,j}(v^e)
\end{eqnarray}

We will now use the $\GL_r$-bundle $\xi$ over $C_{\gen}$,
the isomorphisms (\ref{1}) and (\ref{3}), the numbers
$a_1,\dots,a_r$ and the matrices $H^1(s)$ and $H^2(t)$,
to construct a Gieseker vector bundle of rank $r$ on the curve $C$.

Let $p_1$ and $p_2$ denote the closed point of 
$\Spec(k[[s]])$ and $\Spec(k[[t]])$ respectively.
Let $\V$ be the trivial vector bundle $\Oo^{[1,r]}$
on the disjoint union 
$\Spec(k[[s]])\sqcup\Spec(k[[t]])$ 
(the normalization of $\Spec(k[[s,t]]/(s\cdot t))$),
and let $V$ and $W$ be its fiber at $p_1$ and $p_2$ respectively.
Of course, both $V$ and $W$ are naturally identified with $k^{[1,r]}$.

The numbers $a_1,\dots,a_r$ define a partition
$$
[1,r]=D_1\sqcup D_2\sqcup\dots\sqcup D_m
\quad
$$
characterized by the following properties:
\begin{enumerate}
\item
$D_1$ is the (possibly empty) 
set of all indices $i$ such that $a_i=0$.
\item
For $\nu\geq 2$ the set $D_\nu$ is non-empty.
\item
If $1\leq\nu<\nu'\leq m$, $i\in D_\nu$ and $j\in D_{\nu'}$ 
then $a_i<a_j$.
\item
For all $\nu\in[1,m]$ and $i,j\in D_\nu$ we have $a_i=a_j$.
\end{enumerate}

We define filtrations
\begin{eqnarray*}
0= & F_0(V)\subseteq F_1(V)\subsetneq F_2(V)\subsetneq\dots
                        \subsetneq F_{m-1}(V)
                        \subsetneq F_{m}(V) & =V
\\
0=& F_0(W)\subsetneq F_1(W)\subsetneq F_2(W)\subsetneq\dots
                         \subsetneq F_{m-1}(W)
                         \subseteq F_{m}(W) & =W
\end{eqnarray*}
by setting
$$
F_i(V) := k^{D_{1}\sqcup\dots\sqcup D_i}
\qquad\text{and}\qquad
F_i(W) := k^{D_{m-i+1}\sqcup\dots\sqcup D_m} 
$$
for $i=0,\dots,m$.
For $i=1,\dots,m-1$ let
$$
\varphi_i: F_{m-i}(W)/F_{m-i-1}(W)=k^{D_{i+1}} 
           \ \Isomto\  
           k^{D_{i+1}}=F_{i+1}(V)/F_{i}(V)
$$
be the identity morphism on $k^{D_{i+1}}$ and let 
$\varphib_i$ be the homothety class of $\varphi_i$.
Finally let 
$$
\Phi': F_1(V)/F_{0}(V)= k^{D_1} \ \Isomto\ k^{D_1}=F_m(W)/F_{m-1}(W)
$$
be the identity morphism on $k^{D_1}$.

By \cite{kgl} 9.3 the data 
$$
((F_{\bullet}(V),F_{\bullet}(W)),\varphib_1,\dots,\varphib_{m-1},\Phi')
$$
define a $k$-valued point of $KGL(V,W)$, i.e.
a generalized isomorphism $\Phi$ from $V$ to $W$.

Let $\Ct\to C$ be the normalization of the curve $C$.
By a slight abuse of notation we denote also by $p_1$, $p_2$ the
two points of $\Ct$ which lie above the singular point $p$ of $C$.
Let $\E_{\xi}$ be the rank $r$ vector bundle on 
$C_{\gen}=\Ct\setminus\{p_1,p_2\}$
associated to the principal $\GL_r$-bundle $\xi$.

We use the isomorphism
$$
\xymatrix{
(k((s))\times k((t)))^{[1,r]}
\ar[r]^{(H^1,H^2)}
&
(k((s))\times k((t)))^{[1,r]}
\\
\V\tensor_{k[[s]]\times k[[t]]} (k((s))\times k((t)))
\ar@{=}[u]
&
\E_{\xi}\tensor_{\Oo_{\Ct}}\widehat{K}_p
\ar[u]^{\isomorph}_{\text{(\ref{1}), (\ref{3})}}
}
$$
as a glueing datum to define a vector bundle $\E$
on $\Ct$, whose fibers at the points $p_1$ and $p_2$ are naturally
identified with $V$ and $W$ respectively.
By \ref{GVBD->KGL} the pair $(\E,\Phi)$ induces a 
Gieseker vector bundle datum $(C'\to C,\F,p')$ on $(\Ct,p_1,p_2)$ which
in turn induces a Gieseker vector bundle $(C'\to C, \F)$ on $C$.

For the convenience of the reader I will now describe the Gieseker
vector bundle $(C'\to C, \F)$ explicitely.
Let $R_0:=\Spec(k[[s]])$, $R_m:=\Spec(k[[t]])$.
If $m=1$, we set $R=\Spec(k[[s,t]]/(s\cdot t))$, which is nothing
else but the nodal curve which arrises from $R_0\sqcup R_m$ by
identifying the points $p_1$ and $p_2$.
If $m\geq 2$,
let $R_1,\dots,R_{m-1}$ be $m-1$ copies of the projective line
$\Pp^1$ and let $x_i, y_i$ be two distinct points
in $R_i$. Let $R$ be the nodal curve which arrises from the
union
$$
R_0\sqcup R_1\sqcup\dots\sqcup R_{m-1}\sqcup R_m
$$
by identifying $p_1\in R_0$ and $p_2\in R_m$ with 
$x_1\in R_1$ and $y_{m-1}\in R_{m-1}$ respectively
and by identifying $y_i\in R_i$ with $x_{i+1}\in R_{i+1}$ for $i\in[1, m-2]$:

\vspace{1mm}

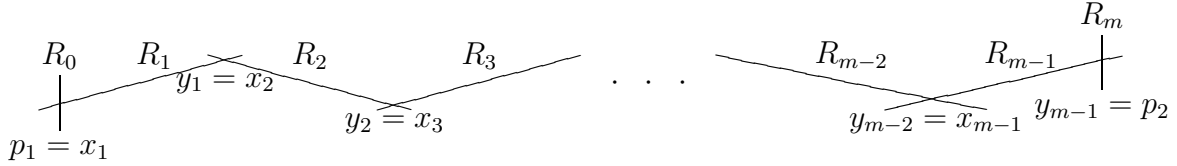
\begin{center}
\parbox{8cm}{
\xy <0mm,-4mm>; <.9mm,-4mm> :
(3,1);(3,-7)**@{-};
(3,-10) *{\text{$p_1=x_1$}};
(3,4) *{\text{$R_0$}};
(0,-4);(30,4)**@{-}; (17,4) *{\text{$R_1$}};
(27.5,0) *{\text{$y_1=x_2$}};
(25,4);(55,-4)**@{-}; (40,4) *{\text{$R_2$}};
(52.5,-6) *{\text{$y_2=x_3$}};
(50,-4);(80,4)**@{-}; (65,4) *{\text{$R_3$}};
(85,0)*{\cdot};
(90,0)*{\cdot};
(95,0)*{\cdot};
(100,4);(140,-4) **@{-}; (120,4) *{\text{$R_{m-2}$}};
(132.5,-6) *{\text{$y_{m-2}=x_{m-1}$}};
(125,-4);(160,4) **@{-}; (145,4) *{\text{$R_{m-1}$}};
(157,-1);(157,7)**@{-};
(157,10) *{\text{$R_m$}};
(157,-4) *{\text{$y_{m-1}=p_2$}};

\endxy
}
\end{center}

\vspace{3mm}

Let $\Oo_{R_i}(1)$ be the defining bundle on $R_i=\Pp^1$
together with isomorphisms
\begin{eqnarray}
\label{triv O(1)} 
\Oo_{R_i}(1)[x_i]\isomto k
\quad\text{and}\quad
\Oo_{R_i}(1)[y_i]\isomto k
\quad.
\end{eqnarray}
We define the rank $r$ vector bundles
$$
E_i:=\Oo_{R_i}^{D_1\sqcup\dots\sqcup D_{i}} \oplus
     \Oo_{R_i}(1)^{D_{i+1}} \oplus
     \Oo_{R_i}^{D_{i+2}\sqcup\dots\sqcup D_{m}}
$$
on $R_i$ together with  the isomorphisms
\begin{eqnarray}
\label{triv E} 
E_i[x_i]\isomto k^{[1,r]}
\quad\text{and}\quad
E_i[y_i]\isomto k^{[1,r]}
\end{eqnarray}
induced by (\ref{triv O(1)}).

The maximal ideal $sk[[s]]$ of $k[[s]]$ is a free module of rank one
and as such defines a line bundle $\Oo_{R_0}(-1)$
on $R_0=\Spec(k[[s]])$. 
We consider this line bundle together with the isomorphism
\begin{eqnarray}
\label{triv O(-1)} 
\Oo_{R_0}(-1)[p_2]\isomto k
\end{eqnarray}
given by $sk[[s]]/s^2k[[s]]\to k$, $s\mapsto 1$.
The generic fiber of $\Oo_{R_0}(-1)$ is identified with $k((s))$ via
the inclusion $sk[[s]]\injto k[[s]]$.
Then we have the rank $r$ vector bundles
$$
E_0:=\Oo_{R_0}^{D_1}\oplus\Oo_{R_0}(-1)^{D_2\sqcup\dots\sqcup D_m}
\qquad\text{and}\qquad
E_m:=\Oo_{R_m}^{[1,r]}
$$
on $R_0$ and $R_m$ together with isomorphisms
\begin{eqnarray}
\label{triv E_0} 
E_0[p_1]\isomto k^{[1,r]}
\quad\text{and}\quad
E_m[p_2]\isomto k^{[1,r]}
\end{eqnarray}
(the first one being induced by (\ref{triv O(-1)})) and isomorphisms
\begin{eqnarray}
\label{triv E_0 gen} 
E_0\tensor_{\Oo_{R_0}}k((s))\isomto k((s))^{[1,r]}
\quad\text{and}\quad
E_m\tensor_{\Oo_{R_m}}k((t))\isomto k((t))^{[1,r]}
\quad.
\end{eqnarray}
The vector bundles $E_0,\dots,E_m$ glue together via the isomorphisms
(\ref{triv E}) and (\ref{triv E_0}) to form a rank $r$ 
vector bundle $E$ on $R$.

Let $C'\to C$ be the modification of $C$ obtained by glueing together
$R$ and $C_{\gen}$ along the isomorphism 
$$
\xymatrix@R=2ex{
\Spec(k((s)))\sqcup\Spec(k((t))) \ar[r]^(.7){(\ref{1})} \ar@{^(->}[d]
&
\Spec(\widehat{K}_p) \ar@{^(->}[d]
\\
R & C_{\gen} 
}
$$
and let $\F$ be the rank $r$ vector bundle on $C'$ obtained by
glueing together $E$ and $\E_{\gen}$ via the isomorphism
$$
\xymatrix{
(k((s))\times k((t)))^{[1,r]}
\ar[r]^{(H^1,H^2)}
&
(k((s))\times k((t)))^{[1,r]}
\\
E\tensor_{\Oo_{R}} (k((s))\times k((t)))
\ar[u]^{\isomorph}_{\text{(\ref{triv E_0 gen})}}
&
\E_{\xi}\tensor_{\Oo_{\Ct}}\widehat{K}_p
\ar[u]^{\isomorph}_{\text{(\ref{1}), (\ref{3})}}
}
$$
It is easy to check that $(C'\to C, \F)$ is indeed a Gieseker
vector bundle on $C$.

\vspace{3mm}
It remains to be shown that the association
$$
(\xi, C\to\Spec(k), \Ac)
\ \mapsto\ 
(C'\to C, \F)
$$
is surjective and
does not depend on the choices 
(\ref{1}), (\ref{2}), (\ref{3}), (\ref{4})
which we made during the construction.
This will be done in the next sections.

\vspace{5mm}
\section{Independence of the isomorphisms (\ref{1}) and (\ref{2})}

Let 
$$
\widehat{\Oo}_p\isomto k[[s,t]]/(s\cdot t)
\eqno(\ref{1}')
$$
be another isomorphism and let
$$
\widehat{\Oo}_q\isomto k[[u,v]]/(u\cdot v)
\eqno(\ref{2}')
$$
be an isomorphism with the required property with respect to 
(\ref{1}').
For the moment we make the following assumption:
$$
\text{
The images of the
two minimal ideals of $\widehat{\Oo}_p$ 
under (\ref{1}) and (\ref{1}')
are the same. 
}
\eqno(*)
$$
Then there are units $\sigma(s),\pi(s)\in k[[s]]^\times$ and 
$\tau(t),\omega(t)\in k[[t]]^\times$ such that the following diagrams commute:
$$
\xymatrix{
k[[s,t]]/(s\cdot t) 
\ar[rr]^{s \mapsto s\sigma(s)}_{t\mapsto t\tau(t)}
& &
k[[s,t]]/(s\cdot t)
\\
&
\text{$\widehat{\Oo}_p$} 
\ar[lu]^{(1)} 
\ar[ru]_{(1')} 
&
}
$$

$$
\xymatrix{
k[[u,v]]/(u\cdot v) 
\ar[rr]^{u \mapsto u\pi(u^e)}_{v\mapsto v\omega(v^e)}
& &
k[[u,v]]/(u\cdot v)
\\
&
\text{$\widehat{\Oo}_p$} 
\ar[lu]^{(2)} 
\ar[ru]_{(2')} 
&
}
$$
Furthermore we have $\pi^e=\sigma$ and $\omega^e=\tau$.

It should be noticed that the $e$-th root of unity $\zeta$ is independent
of whether we choose (\ref{1}) or (\ref{1}'),
since it is the eigenvalue of $\gamma$ operating on the tangent
space of one of the branches of $\Spec(\widehat{\Oo}_q)$
and by assumption $(*)$ both the isomorphisms 
(\ref{1}) and (\ref{1}') map  that branch 
$\Spec(\widehat{\Oo}_q)$ to the same branch $\{v=0\}$ of
$\Spec(k[[u,v]]/(u\cdot v))$.
Therefore the elements $\alpha_1,\dots,\alpha_r\in\Z/e\Z$ and 
the numbers $a_1,\dots,a_r$ are independent
of whether we choose (\ref{1}) or (\ref{1}').

Let $(\tF_1(u),\tF_2(v))$ be the image of $F$ under the isomorphism
$
\GL_r(\widehat{K}_q)\isomto\GL_r(k[[u]])\times\GL_r(k[[v]])
$
induced by (\ref{2}'). Then we have $\tF^1(u)=F^1(u\pi(u^e))$ and
$\tF^2(v)=F^2(v\omega(v^e))$ and it follows that
\begin{eqnarray*}
\tF^1_{i,j}(u) &=& u^{a_j}\cdot\tH^1_{i,j}(u^e)
\quad,
\\
\tF^2_{i,j}(v) &=& v^{-a_j}\cdot\tH^2_{i,j}(v^e)
\quad, 
\end{eqnarray*}
where
\begin{eqnarray*}
\tH^1_{i,j}(s) &=& \pi^{a_j}H^1_{i,j}(s\sigma)
\quad,
\\
\tH^2_{i,j}(t) &=& \omega^{-a_j}H^2_{i,j}(t\tau)
\quad.
\end{eqnarray*}
Therefore the following diagram commutes:
$$
\xymatrix{
(k((s))\times k((t)))^{[1,r]}
\ar[rr]^{(H^1,H^2)} 
\ar[dd]
& &
(k((s))\times k((t)))^{[1,r]}
\ar[dd]_{s\mapsto s\sigma}^{t\mapsto t\tau}
\\
&
\E_{\xi}\tensor_{\Oo_{\Ct}}\widehat{K}_p
\ar[ur]^{\isomorph}_{\text{(\ref{1}), (\ref{3})}}
\ar[dr]_{\isomorph}^{\text{(\ref{1}'), (\ref{3})}}
&
\\
(k((s))\times k((t)))^{[1,r]}
\ar[rr]^{\text{$(\tH^1,\tH^2)$}}
& &
(k((s))\times k((t)))^{[1,r]}
}
$$
where the left vertical arrow maps an element $(x(s),y(t))$
to the element 
$$
(\diag(\pi^{a_1},\dots,\pi^{a_r}) x(s\sigma),
 \diag(\omega^{-a_1},\dots,\omega^{-a_r}) y(t\tau))
\quad.
$$
Let $\tilde{\E}$ be the vector bundle on $\Ct$ obtained by the glueing datum
$(\tH^1,\tH^2)$. Then the above diagram shows that there is an
isomorphism  $\E\isomto\tilde{\E}$ which induces the isomorphisms
$$
\xymatrix@R=1ex{
\E[p_1]=k^{[1,r]} 
\ar[rrrr]^{\diag(\pi(0)^{-a_1},\dots,\pi(0)^{-a_r})}
& & & &
k^{[1,r]}=\tilde{\E}[p_1]
\\
\E[p_2]=k^{[1,r]} 
\ar[rrrr]^{\diag(\omega(0)^{a_1},\dots,\omega(0)^{a_r})}
& & & &
k^{[1,r]}=\tilde{\E}[p_2]
}
$$
between the fibers at $p_1$ and $p_2$ respectively.
Thus it maps the generalized isomorphism 
$
\Phi\in\KGL(k^{[1,r]},k^{[1,r]})=\KGL(\E[p_1],\E[p_2])
$
to the generalized isomorphism
$
\Phi\in\KGL(k^{[1,r]},k^{[1,r]})=\KGL(\tilde{\E}[p_1],\tilde{\E}[p_2])
$.
This shows that the pairs $(\E,\Phi)$
and $(\tilde{\E},\Phi)$ are isomorphic. Consequently
this is also true for the associated Gieseker vector bundles.

To get rid of the assumption $(*)$ we investigate now what happens
if we change the isomorphisms 
(\ref{1}), (\ref{2})  
by composing them with the automorphisms
$$
\xymatrix@R=1ex{
k[[s,t]]/(s\cdot t) \ar[r]^{s\mapsto t}_{t\mapsto s}
&
k[[s,t]]/(s\cdot t)
\\
k[[u,v]]/(u\cdot v) \ar[r]^{u\mapsto v}_{v\mapsto u}
&
k[[u,v]]/(u\cdot v)
}
$$
respectively.

This means that $\zeta^{-1}$ takes the role of $\zeta$ and consequently
the set $\{\alpha_1,\dots,\alpha_r\}\subseteq \Z/e\Z$ from 
\ref{alpha} is replaced by the set $\{-\alpha_1,\dots,-\alpha_r\}$.
It follows that in (\ref{8}) we would choose integers
$
\tilde{a}_1,\dots,\tilde{a}_r
$
instead of
$
a_1,\dots,a_r
$,
where
$$
\tilde{a}_i=
\left\{
\begin{array}{ll}
a_i \quad & \text{for $i\in[1,i_1]=D_1$} \\
e-a_{r+i_1+1-i} \quad & \text{for $i\in[i_1+1,r]$} 
\end{array}
\right.
$$
Then the matrix $F$ is replaced by the matrix 
$\tilde{F}=F\cdot\Lambda$, where 
$$
\Lambda =
\left[
\vcenter{
\hbox{
\xymatrix@C.5ex@R.5ex{
\eins_{i_1} & \ar@{-}'[dddd] &   & 0 &    \\
\ar@{-}'[rrrr]          &   &   &   &    \\ 
            &   &   &   &  1 \ar@{.}[ddll] \\
     0      &   &   &   &    \\
            &   & 1 &   &   
}
}}
\right]
$$
is the permutation matrix
belonging to the permutation $\lambda\in S_r$, where
$$
\lambda(i)=
\left\{
\begin{array}{ll}
i \quad & \text{for $i\in[1,i_1]$} \\
r+i_1+1-i \quad & \text{for $i\in[i_1+1,r]$}
\end{array}
\right.
\quad,
$$
and the matrices $H^1(s)$ and $H^2(t)$ are replaced by the matrices
$\tilde{H}^1(s)$ and $\tilde{H}^2(t)$ respectively, where
$$
\tilde{H}^1(s) = H^2(s)\cdot\Lambda\cdot
\left[
\vcenter{
\hbox{
\xymatrix@C.5ex@R.5ex{
\eins_{i_1} & 0 \\
0 & s^{-1}\eins_{r-i_1}
}
}}
\right]
\qquad\text{and}\qquad
\tilde{H}^2(t) = H^1(t)\cdot\Lambda\cdot
\left[
\vcenter{
\hbox{
\xymatrix@C.5ex@R.5ex{
\eins_{i_1} & 0 \\
0 & t\eins_{r-i_1}
}
}}
\right]
\quad.
$$
The numbers $\tilde{a}_1,\dots,\tilde{a}_r$ define the partition
$$
[1,r]=\tilde{D}_1\sqcup\tilde{D}_2\sqcup\dots\sqcup\tilde{D}_m
$$
where $\tilde{D}_1=D_1$ and $\tilde{D}_i=\lambda(D_{m+2-i})$
for $i\in[2,m]$.
Let 
$(\tilde{R}=
  \tilde{R}_0\cup\dots\cup\tilde{R}_m,\tilde{E})$ 
be the nodal curve associated to this partition,
together with isomorphisms 
$\tilde{E}\tensor_{\Oo_{\tilde{R}}}k((s))\isomto k((s))^{[1,r]}$
and
$\tilde{E}\tensor_{\Oo_{\tilde{R}}}k((t))\isomto k((t))^{[1,r]}$
as in (\ref{triv E_0 gen}).

Now one checks easily that there is an isomorphism
$$
\rho:(R,E)\isomto(\tilde{R},\tilde{E})
$$
which sends the component $R_i$ to $\tilde{R}_{m-i}$ ($i=0,\dots,m$),
such that the following diagram commutes:
$$
\xymatrix{
E\tensor_{\Oo_R}(k((s))\times k((t)))
\ar[r]^{\isomorph}
\ar[d]^{\rho}
&
(k((s))\times k((t)))^{[1,r]}
\ar[d]_{\rho'}
\\
\text{$\tilde{E}\tensor_{\Oo_{\tilde{R}}}(k((s))\times k((t)))$}
\ar[r]^{\isomorph}
&
(k((s))\times k((t)))^{[1,r]}
}
$$
where the morphism $\rho'$ is given by
$$
(x(s),y(t))\mapsto
\left(
\Lambda\cdot
\left[
\vcenter{
\hbox{
\xymatrix@C.5ex@R.5ex{
\eins_{i_1} & 0 \\
0 & s^{-1}\eins_{r-i_1}
}
}}
\right]\cdot
y(s)
\ ,\ 
\Lambda\cdot
\left[
\vcenter{
\hbox{
\xymatrix@C.5ex@R.5ex{
\eins_{i_1} & 0 \\
0 & t\eins_{r-i_1}
}
}}
\right]\cdot
x(t)
\right)
$$
From the commutativity of the diagram
$$
\xymatrix{
(k((s))\times k((t)))^{[1,r]}
\ar[rr]^{(H^1(s),H^2(t))}
\ar[d]^{\rho'}
& &
(k((s))\times k((t)))^{[1,r]}
\ar[d]^{s\mapsto t}_{t\mapsto s}
&
\E_{\xi}\tensor_{\Oo_{\tilde{C}}}\widehat{K}_p
\ar[l]_(.4){(1),(3)}
\ar@{=}[d]
\\
(k((s))\times k((t)))^{[1,r]}
\ar[rr]^{(\tilde{H}^1(s),\tilde{H}^2(t))}
& &
(k((s))\times k((t)))^{[1,r]}
&
\E_{\xi}\tensor_{\Oo_{\tilde{C}}}\widehat{K}_p
\ar[l]_(.4){(1'),(3)}
}
$$
it finally follows that the Gieseker vector bundle
$(\tilde{C}'\to C,\tilde{\F})$
constructed from the data $\xi$, (\ref{1}'),
(\ref{3}), $(\tilde{a}_1,\dots,\tilde{a}_r)$, 
$\tilde{H}^1(s)$, $\tilde{H}^2(t)$
is isomorphic to the Gieseker vector bundle 
$(C'\to C,\F)$ constructed from the data $\xi$, (\ref{1}),
(\ref{3}), $(a_1,\dots,a_r)$, $H^1(s)$, $H^2(t)$.

\vspace{5mm}
\section{Independence of the isomorphisms (\ref{3}) and (\ref{4})}

Independence of (\ref{3}) is immediate, since if we change it by
an automorphism of $\GL_r\times\Spec(\widehat{K}_p)$
(which can be written as an element in $\GL_r(k((s)))\times\GL_r(k((t)))$), 
then $(H^1(s),H^2(t))$ is changed by that same matrix.

Two isomorphisms (\ref{4}) differ by a matrix
$
A=(A_{i,j}) \in \GL_r(\widehat{\Oo}_q)
$
such that
\begin{eqnarray}
\label{cond on A}
A = \diag(\zeta^{-\alpha_1},\dots,\zeta^{-\alpha_r})\cdot
    \gamma(A)\cdot
    \diag(\zeta^{\alpha_1},\dots,\zeta^{\alpha_r})
\quad.
\end{eqnarray}
After identifying $\widehat{\Oo}_q$ with the ring 
$k[[u,v]]/(u\cdot v)$ via the isomorphism (\ref{1}), we can write 
$$
A=A^0 + u\cdot A^1(u) + v\cdot A^2(v)
$$
with uniquely determined matrices
$A^0\in\GL_r(k)$, 
$A^1(u)\in M(r\times r, k[[u]])$ and 
$A^2(v)\in M(r\times r, k[[v]])$.
Condition (\ref{cond on A}) implies that $A^0$ is a block matrix 
of the form
\begin{equation}
\label{A^0}
A^0 = 
\left[
\vcenter{
\hbox{
\xymatrix@C.5ex@R.5ex{
A^0_1 \ar@{.}[dr] & 0 \\
0 & A^0_m
}
}}
\right]
\end{equation}
where the $A^0_i$ are blocks of size $\card D_i$ for $i=1,\dots,m$.
Condition (\ref{cond on A}) implies furthermore that there
are matrices 
$B^1(s)=(B^1_{i,j}(s))\in\GL_r(k[[s]])$ and 
$B^2(t)=(B^2_{i,j}(t))\in\GL_r(k[[t]])$
such that
\begin{eqnarray*}
A^1(u) &=& u^{-1}\diag(u^{a_1},\dots,u^{a_r}) 
           \cdot B^1(u^e)\cdot
           \diag(u^{-a_1},\dots,u^{-a_r}) 
\quad,
\\
A^2(v) &=& v^{-1}\diag(v^{-a_1},\dots,v^{-a_r}) 
           \cdot B^2(v^e)\cdot
           \diag(v^{a_1},\dots,v^{a_r}) 
\end{eqnarray*}
and such that 
\begin{equation}
\label{B(0)}
\begin{array}{ll}
B^1_{i,j}(0) &= 0 
\quad\text{for $a_i-a_j\leq 0$}
\quad,
\\
B^2_{i,j}(0) &= 0 
\quad\text{for $a_j-a_i\leq 0$}
\quad.
\end{array}
\end{equation}

The change of (\ref{4}) by the matrix $A$ means that we have
to replace $F$ by the matrix 
$$
\tilde{F}=F\cdot A
$$
and that consequently we have to replace
the matrices $H^1(s)$ and $H^2(t)$ by the matrices
\begin{eqnarray*}
\tilde{H}^1(s) &=& H^1(s)\cdot(A^0 + B^1(s)) 
\qquad\text{and}
\\
\tilde{H}^2(t) &=& H^2(t)\cdot(A^0 + B^2(t))
\end{eqnarray*}
respectively.

The pair of matrices 
$(A^0 + B^1(s)$, $A^0 + B^2(t))$ 
defines an automorphism of $\V$ which induces the automorphisms
$A^0 + B^1(0)$ and $A^0 + B^2(0)$ on the special fibers
$V$ and $W$ respectively.
From (\ref{A^0}) and (\ref{B(0)}) it follows that the induced automorphism
of $\KGL(V,W)$ maps the generalized isomorphism $\Phi$ to itself.

It follows that the pair $(\tilde{\E},\tilde{\Phi})$ obtained by the 
glueing datum $(\tilde{H}^1, \tilde{H}^2)$ is isomorphic to the pair
$(\E,\Phi)$ obtained by the glueing datum $(H^1, H^2)$. Therefore
also the induced Gieseker vector bundles are isomorphic.

\vspace{5mm}
\section{Surjectivity}

Let $(C'\to C,\F)$ be a Gieseker vector bundle on $C$.
By definition, $C'$ is either isomorphic to $C$, or it is
the union of the normalization $\Ct$ of $C$ and a chain
$R$ of projective lines which intersects $\Ct$ in the two
points $p_1$ and $p_2$ lying above the singularity $p\in C$.
In the first case we let $p':=p$, in the second case
we let $p'=p_2$.
Then the tripel
$$
(\Ct'\to\Ct, \Ft', p')
$$
is a Gieseker vector bundle datum in the sense of \ref{GVBD}.
By \ref{GVBD->KGL} such a datum induces a vector bundle
$\E$ on the curve $\Ct$ together with a generalized isomorphism 
$\Phi$ from $V:=\E[p_1]$ to $W:=\E[p_2]$.

More precisely, $\Phi$ is a $k$-valued point of $\KGL(V,W)$ which
lies in the stratum $\bO_{I,J}$ for some $I\subseteq[0,r-1]$ and $J=\emptyset$.
As we have recalled in \S \ref{review}, such a point is given by a tupel
$$
((F_{\bullet}(V),F_{\bullet}(W)), \varphib_1,\dots,\varphib_{m-1},\Phi')
\quad,
$$
where $m:=|I|+1$, 
\begin{eqnarray*}
0= & F_0(V)\subseteq F_1(V)\subsetneq F_2(V)\subsetneq\dots
                        \subsetneq F_{m-1}(V)
                        \subsetneq F_{m}(V) & =V
\\
0=& F_0(W)\subsetneq F_1(W)\subsetneq F_2(W)\subsetneq\dots
                         \subsetneq F_{m-1}(W)
                         \subseteq F_{m}(W) & =W
\end{eqnarray*}
are flags in $V$ and $W$ respectively, 
$\varphib_i$ is the homothety class of an isomorphism
$$
\varphi_{\nu}: F_{m-\nu}(W)/F_{m-\nu-1}(W)\isomto F_{\nu+1}(V)/F_{\nu}(V)
$$
and $\Phi'$ denotes an isomorphism 
$F_1(V)/F_0(V)\isomto F_m(W)/F_{m-1}(W)$.

There is a basis $v_1,\dots,v_r$ of $V$ and $w_1,\dots,w_r$ of $W$
and a partition
$$
[1,r]=D_1\sqcup D_2\sqcup \dots \sqcup D_m
$$
with the property that
\begin{enumerate}
\item
$i\in D_{\nu}$, $j\in D_{\nu'}$ and $i<j$ implies $\nu\leq \nu'$,
\item
$F_{\nu}(V)$ is generated by 
$\{v_i\ |\ i\in D_1\sqcup\dots\sqcup D_{\nu}\}$
and
$F_{\nu}(W)$ is generated by 
$\{w_i\ |\ i\in D_{m-\nu+1}\sqcup\dots\sqcup D_{m}\}$,
\item
for $i\in D_{\nu+1}$
the isomorphism $\varphi_{\nu}$ sends the residue class of $w_i$
mod $F_{m-\nu-1}(W)$
to the residue class of $v_i$ mod $F_{\nu}(V)$,
\item
for $i\in D_1$ the isomorphism $\Phi'$ sends 
$v_i$  to the residue class of $w_i$ mod $F_{m-1}(W)$.
\end{enumerate}

For $i=1,2$ 
we choose an isomorphism 
\begin{equation}
\label{E at p_i}
\E\tensor_{\Oo_{\Ct}}\Oh_{p_i}\isomto \Oh_{p_i}^r
\end{equation}
which induces the isomorphism 
$V\to k^r$, $v_i\to e_i$ 
and
$W\to k^r$, $w_i\to e_i$ 
from the fibres at $p_1$ and $p_2$ respectively,
where $e_1,\dots,e_r$ is the canonical basis of $k^r$.
Let $\xi$ be the prinicpal $\GL_r$-bundle on $C_{\gen}$ of local frames of 
the restriction of the vector bundle $\E$ to
$C_{\gen}=\Ct\setminus\{p_1,p_2\}$.
Then the isomorphisms \ref{E at p_i} induce the isomorphism
\begin{equation}
\label{xi at K_p}
\xi\times_{C_{\gen}}\Spec(\Kh_p)
\isomto
\GL_r\times\Spec(\Kh_p)
\end{equation}

\begin{lemma}
There is a morphism $f:U\to C$, an integer  
$e\geq m$ prime to the characteristic of $k$ 
and an operation of $\Gamma:=\Z/e\Z$
on $U$ such that
\begin{enumerate}
\item
$\Gamma$ leaves $f$ invariant and the induced morphism 
$U/\Gamma\to C$ is etale,
\item
$U$ has exactly one singular point $q$ and $f^{-1}(p)=\{q\}$,
\item
the action of $\Gamma$ on $f^{-1}(C_{\gen})$ is free.
\end{enumerate}
\end{lemma}

\begin{proof}
Since $p$ is an ordinary double point of $C$, there exists a diagram
of pointed schemes and \'etale morphisms as follows:
$$
\xymatrix{
(C,p) 
&
(U_0,q_0) 
\ar[l]_{\text{\'etale}}
\ar[r]^(.24){\text{\'etale}}
&
(V_0,y_0):=((\Spec(k[s,t]/(s\cdot t)),(s,t))
}
\quad.
$$
After removing from $U_0$
the points $\neq q_0$ in the fiber of $U_0\to C$
we may assume that $q_0$ is the only point lying above $p$.
Choose $e\in\Z$ prime to $\chara(k)$ with $e\geq m$. 
Let 
$$(V,y):=(\Spec(k[u,v]/(u\cdot v)),(u,v))$$
and let $(V,y)\to(V_0,y_0)$ be defined by $s\mapsto u^e$, $t\mapsto v^e$.
Let $\gamma$ be a generator of $\Gamma=\Z/e\Z$ and let $\zeta\in k$ be a 
primitive $e$-th root of unity. We define an action of $\Gamma$ on $(V,y)$
by letting $\gamma(u)=\zeta u$ and $\gamma(v)=\zeta^{-1} v$.
Now we set
$$
U:=U_0\times_{V_0}V
$$
and let $f:U\to C$ be the composition $U\to U_0\to C$.
From $V$ the scheme $U$ inherits an action of the group $\Gamma$.
Since $V/\Gamma=V_0$ and $U_0\to V_0$ is flat we have 
$U/\Gamma=U_0$ which by construction is \'etale over $C$.
The only point in the fibre of $f$ over $p$ is the point
$q=(q_0,y)\in U$. Since $U_0\to V_0$ is \'etale
and $V\to V_0$ is smooth outside the point $y$,
it follows that the fibre product $U=U_0\times_{V_0}V$ is regular
outside $q$.
Furthermore, since the action of $\Gamma$ on $V\setminus\{y\}$ is
free the same holds for the action of $\Gamma$ on $U\setminus \{q\}$.
\end{proof}

In what follows we will construct a chart $(U,\eta,\Gamma)$ for $\xi$
where $U\to C$ and $\Gamma$ are chosen as in the lemma and the 
$\GL_r$-bundle $\eta$ with $\Gamma$-operation 
is glued together from an object $\eta_{\gen}$
over $U_{\gen}$ and an object $\etah_q$ over
the completion of $U$ at the singular point $q$.

To fix notation, 
let $\Oh_q$ be the completion of the local ring $\Oo_{U,q}$ and
let $\gamma$ be a generator of $\Gamma$.
There exists an isomorphism 
\begin{equation}
\Oh_q\isomto k[[u,v]]/(u\cdot v)
\end{equation}
and a primitive $e$-th root of unity $\zeta$ such that
the automorphism $\gamma:\Oh_q\isomto\Oh_q$ translates
into the automorphism $u\mapsto \zeta u$, $v\mapsto \zeta^{-1} v$
of $k[[u,v]]/(u\cdot v)$ (cf. \cite{ACV}, 2.1.2).

Let $a_i\in[0,e-1]$ ($i\in[1,r]$) be chosen such that:
\begin{eqnarray*}
a_i=0    &\text{for}& i\in D_1, \\
a_i<a_j  &\text{for}& i\in D_{\nu},\ j\in D_{\nu'},\ \nu<\nu', \\
a_i=a_j  &\text{for}& i,j\in D_{\nu},\ \nu \in[1,m].
\end{eqnarray*}
Let $\etah_q:=\GL_r\times\Spec\Oh_q$ together with the $\Gamma$-operation
defined by 
$$
\diag(\zeta^{a_1},\dots,\zeta^{a_r})\times\gamma:
\GL_r\times\Spec\Oh_q\isomto\GL_r\times\Spec\Oh_q
\quad.
$$

Let $\eta_{\gen}:=\xi\times_{C_{\gen}}U_{\gen}$ together with the
$\Gamma$-operation given by 
$$
\id\times\gamma:\xi\times_{C_{\gen}}U_{\gen}\isomto
                \xi\times_{C_{\gen}}U_{\gen}
\quad.
$$

Now we glue together $\etah_q$ and $\eta_{\gen}$ along
$\Spec\Kh_q\isomorph k((u))\times k((v))$ via the
isomorphism
$$
\xymatrix{
\text{$\etah_q\times_{\Oh_q}\Spec(\Kh_q)$}
\ar[r]^{(\ref{xi at K_p})}
&
\text{$\GL_r\times\Spec(\Kh_q)$}
\ar[r]^(0.35){F^1\times F^2}
&
\text{$\GL_r\times\Spec(\Kh_q)$}
=
\text{$\eta_{\gen}\times_{U_{\gen}}\Spec(\Kh_q)$}
}
$$
where
$$
F^1=\diag(u^{a_1},\dots,u^{a_r})
\quad\text{and}\quad
F^2=\diag(v^{-a_1},\dots,v^{-a_r})
\quad.
$$
This gives a principal $\GL_r$-bundle $\eta$ on $U$.
From the commutativity of the diagram
$$
\xymatrix{
\text{$\etah_q\times_{\Oh_q}\Spec(\Kh_q)$}
\ar[r]^{\isomorph}
\ar[d]_{\diag(\zeta^{a_1},\dots,\zeta^{a_r})\times\gamma}
&
\text{$\eta_{\gen}\times_{U_{\gen}}\Spec(\Kh_q)$}
\ar[d]^{\id\times\gamma}
\\
\text{$\etah_q\times_{\Oh_q}\Spec(\Kh_q)$}
\ar[r]^{\isomorph}
&
\text{$\eta_{\gen}\times_{U_{\gen}}\Spec(\Kh_q)$}
}
$$
it follows that the $\Gamma$-operation on $\etah_q$ and $\eta_{\gen}$
induces a $\Gamma$-operation on $\eta$.
It is clear from the construction that the triple $(U,\eta,\Gamma)$
forms a chart for $\xi$.

There is a chart $(U_1,\eta_1,\Gamma_1)$ for $\xi$,
where $U_1:=C_{\gen}$, $\eta:=\xi$,
$\Gamma:=(1)$. 
This chart together with the chart $(U,\eta,\Gamma)$ make up
a balanced atlas $\Ac$ for $\xi$.
It is clear by construction that the twisted $G$-bundle
$(\xi,C\to\Spec(k),\Ac)$ is mapped to the Gieseker
vector bundle $(C'\to C,\F)$.

\section{Further directions}

The relationship between twisted $GL_r$-bundles and
Gieseker vector bundles should be further investigated since
it might lead to a clue what the right notion of stable maps 
to the classifying stack of a reductive group are.
The next step would be to try to extend the mapping given in
\ref{thm2} so that it works for families.

For example let $A:=\C[[t]]$, let $S:=\Spec A$ 
and let $C\to S$ be a 
stable curve over $S$. Let
$$
\left(
\vcenter{
\xymatrix@R=1ex{
C' \ar[rr] \ar[dr] & &
C \ar[dl] \\
& S &
}},
\F
\right)
$$
be a Gieseker vector bundle of rank $r$ on $C$.

Assume in paticular that the generic fiber of $C\to S$ is smooth
and that its special fiber is irreducible with one double point $p$.
Then it can be shown that there is a twisted $GL_r$-bundle 
$(\xi,C\to S,\Ac)$ such that if we apply the mappings from theorems 
\ref{thm1} and \ref{thm2} to the isomorphism class of the 
generic and special fiber of $(\xi,C\to S,\Ac)$, 
then we obtain the generic and special fiber of the Gieseker vector
bundle
$(C'\to C,\F)$ respectively.
Indeed, in the neighbourhood of $p$ one may chose a chart 
$(U,\eta,\Gamma)$ for $\xi$, where $U\to C$ \'etale locally looks
like 
$$
\xymatrix@R=1ex{
\Spec A[u,v]/(uv-t)\ar[r] &
\Spec A[x,y]/(xy-t^e) \\
u^e &
x \ar@{|->}[l] \\
v^e &
y \ar@{|->}[l]
}
$$

On the other hand, assume $C=C_0\times S$, where $C_0$ is an irreducible
curve with one ordinary double point $p$, and assume that 
$C'\to C$ induces an isomorphism of the generic fibers and the
morphism $C_1\to C_0$ on the special fibers.
For simplicity let us assume furthermore that the rank $r$ of
the Gieseker bundle is one.
In this situation it would be interesting to
know, whether there is a twisted $\GL_1$-bundle
$(\xi,C\to S,\Ac)$ such that the map of theorem \ref{thm2}
maps the generic and the special fiber of $(\xi,C\to S,\Ac)$
to the generic and special fiber of $(C'\to C,\F)$ respectively.


\end{document}